\documentclass{gtmon_a}
\pdfoutput=1
\usepackage{pinlabel}

\proceedingstitle{The Zieschang Gedenkschrift}
\conferencestart{5 September 2007}
\conferenceend{8 September 2007}
\conferencename{Conference in honour of Heiner Zieschang}
\conferencelocation{Toulouse, France}

\editor{Michel Boileau}
\givenname{Michel}
\surname{Boileau}

\editor{Martin Scharlemann}
\givenname{Martin}
\surname{Scharlemann}

\editor{Richard Weidmann}
\givenname{Richard}
\surname{Weidmann}

\dedicatory{Dedicated to the memory of Heiner Zieschang, first to notice that genus two handlebodies could be interesting}

\title{Refilling meridians in a genus 2 handlebody complement}   

\author{Martin Scharlemann}
\givenname{Martin}
\surname{Scharlemann}
\address{Mathematics Department\\University of California\\\newline
Santa Barbara, CA USA}
\email{mgscharl@math.ucsb.edu}
\urladdr{}

\volumenumber{14}
\issuenumber{}
\publicationyear{2008}
\papernumber{20}
\startpage{451}
\endpage{475}

\doi{}
\MR{}
\Zbl{}


\received{30 March 2006}
\revised{26 April 2007}
\accepted{26 April 2007}
\published{29 April 2008}
\publishedonline{29 April 2008}
\proposed{}
\seconded{}
\corresponding{}
\version{}

\arxivreference{math/0603705}

\makeatletter
\def\cnewtheorem#1[#2]#3{\newtheorem{#1}{#3}[section]
\expandafter\let\csname c@#1\endcsname\c@thm}
\makeatother

\let\xysavmatrix\xymatrix
\def\xymatrix{\disablesubscriptcorrection\xysavmatrix}
\AtBeginDocument{\let\bar\wbar\let\tilde\wtilde\let\hat\what}

\newtheorem{thm}{Theorem}[section]  
\cnewtheorem{cor}[thm]{Corollary}
\cnewtheorem{lemma}[thm]{Lemma} 
\cnewtheorem{prop}[thm]{Proposition} 
\theoremstyle{remark}
\newtheorem{conj}{Conjecture} 
\cnewtheorem{defin}[thm]{Definition} 
\cnewtheorem{quest}[thm]{Question}
\cnewtheorem{example}[thm]{Example}

\newcommand{\aaa}{\alpha}

\newcommand{\bbb}{\beta}
\newcommand{\sss}{\sigma}  
\newcommand{\Sss}{\Sigma} 
\newcommand{\Uuu}{\Upsilon}

\newcommand{\rrr}{\rho}

\newcommand{\Maa}{M[\aaa]}
\newcommand{\Mbb}{M[\bbb]}
 
\newcommand{\bdd}{\partial}

\newcommand{\Rpp}{\mathbb{RP}^3}
\newcommand{\inter}{\mathrm{int}}

\begin{document}

\begin{asciiabstract}
Suppose a genus two handlebody is removed from a 3-manifold M and then
a single meridian of the handlebody is restored.  The result is a knot
or link complement in M and it is natural to ask whether geometric
properties of the link complement say something about the meridian
that was restored.  Here we consider what the relation must be between
two not necessarily disjoint meridians so that restoring each of them
gives a trivial knot or a split link.
\end{asciiabstract}

\begin{htmlabstract}
Suppose a genus two handlebody is removed from a 3&ndash;manifold M and
then a single meridian of the handlebody is restored.  The result is a
knot or link complement in M and it is natural to ask whether
geometric properties of the link complement say something about the
meridian that was restored.  Here we consider what the relation must
be between two not necessarily disjoint meridians so that restoring
each of them gives a trivial knot or a split link.
\end{htmlabstract}

\begin{abstract} 
Suppose a genus two handlebody is removed from a $3$--manifold $M$ and
then a single meridian of the handlebody is restored.  The result is a
knot or link complement in $M$ and it is natural to ask whether
geometric properties of the link complement say something about the
meridian that was restored.  Here we consider what the relation must
be between two not necessarily disjoint meridians so that restoring
each of them gives a trivial knot or a split link.
\end{abstract}

\maketitle

\section{Background}

For a knot or link in a $3$--manifold, here are some natural geometric
questions that arise, in roughly ascending order of geometric
sophistication: is the knot the unknot? is the link split?  is the
link or knot a connected sum?  are there companion tori?  beyond
connected sums, are there essential annuli in the link complement?
beyond connected sums, are there essential meridional planar surfaces?
One well-established context for such questions is that of Dehn
surgery (cf Gordon \cite{Go}) where one imagines filling in the knot
or link complement with solid tori via different meridian slopes and
then asks under what conditions two of the fillings have geometric
features such as those listed above.

Another natural context is this:  Suppose $W$ is a genus $2$ handlebody embedded in a compact orientable $3$--manifold $M$.  Suppose $\aaa, \bbb$ are not necessarily disjoint essential properly embedded disks in $W$ (called therefore meridian disks).  Then $W - \eta(\aaa)$ (resp $W -  \eta(\bbb)$) is a regular neighborhood of a knot or link $L[\aaa]$ (resp $L[\bbb]$) in $M$.  Under what circumstances do $L[\aaa]$ and $L[\bbb]$ have geometric features like those outlined above?  At the most primitive level (and so presumably the easiest level) one can ask when both $L[\aaa]$ and $L[\bbb]$ are split (if a link) or trivial (if a knot).  Put another way, suppose $M[\aaa], M[\bbb]$ are the manifolds obtained from $M - W$ by restoring neighborhoods of the meridians $\aaa$ and $\bbb$.   Under what circumstances are both $M[\aaa]$ and $M[\bbb]$ reducible and/or $\bdd$--reducible?  (In the absence of Lens space or $S^1 \times S^2$ summands in a closed $3$--manifold $M$, a $\bdd$--reducing disk for a knot complement is equivalent to an unknotting disk for the knot.)  

Not only is this a natural question, the solution to it in specific cases has been of significant interest in knot theory over the past two decades.  Here is a probably partial list of such results, for $M = S^3$:

\begin{itemize}

\item Scharlemann \cite{Sc3} (see also Gabai \cite{Ga}) considers the
case in which one meridian disk is separating, the other is
non-separating, and the two meridians intersect in a single arc.

\item Bleiler and Scharlemann \cite{BS2,BS1} (see also Scharlemann and
Thompson \cite{ST2}) consider the case (among others) in which both
meridian disks are non-separating and the two intersect in a single
arc.

\item Eudave-Mu\~noz \cite{EM3} extends these earlier results and
unifies them within the more general setting in which there are
non-isotopic non-separating meridian disks $\mu, \mu'$ for $W$ which
are disjoint from both $\aaa$ and $\bbb$.

\end{itemize}

If we extend the question to whether one of the links is a connected
sum, the literature becomes even more extensive, including
Eudave-Mu\~noz \cite{EM3,EM2,EM1}, Scharlemann \cite{Sc2}, and
Scharlemann and Thompson \cite{ST1}.

We briefly describe a typical conclusion in the arguments above.  First some terminology:  Say that the handlebody $W \subset M$ is {\em unknotted} if it is isotopic to the regular neighborhood of a figure 8 graph $\Sss$ that lies on a sphere in $M$.  If $W$ is unknotted and $\Sss$ is such a spine and if $\mu$ and $ \mu'$ are the pair of meridian disks for $W$ that are dual to the two edges of $\Sss$, then $\mu$ and $ \mu'$ are called an {\em unknotting pair of meridians for $W$}.  Put another way, the meridians $\mu$ and $\mu'$ for $W$ are called unknotting meridians, and $W$ is said to be unknotted, if there are properly embedded disks $\lambda, \lambda' \subset M - \mathrm{interior}(W)$ such that $|\mu \cap \lambda | = |\mu' \cap \lambda' | = 1$ and $\mu \cap \lambda'  = \mu' \cap \lambda  = \emptyset$.  What is typically proven (most generally by Eudave-Mu\~noz in \cite{EM3}) is this:  We are given specific conditions on the filling meridians $\aaa$ and $\bbb$, including that they are both disjoint from the same pair of non-isotopic meridians $\mu, \mu'$ for $W$.  We suppose further that  the manifolds $M[\aaa]$ and $M[\bbb]$ are both either reducible or $\bdd$--reducible, whereas $M - W$ is irreducible. The conclusion is that $W$ is unknotted in $M$, and the meridians $\mu$ and $\mu'$ are an unknotting pair of meridians for $W$.

Put in this way, one wonders if the various conditions on $\aaa$ and $\bbb$ can be dropped to give a global theorem on the unknottedness of $W$.  For example: 

\medskip
{\bf Naive Conjecture}\qua Suppose $\aaa$ and $\bbb$ are meridian disks
of a genus two handlebody $W \subset M$.  Suppose further that
$\overline{M - W}$ is irreducible and $\bdd M$ is incompressible in
$\overline{M - W}$. Then either $W$ is unknotted or at least one of
$M[\aaa]$, $M[\bbb]$ is both irreducible and $\bdd$--irreducible.

\begin{figure}[ht!]
\labellist\small
\pinlabel $\alpha$ [r] at 21 49
\pinlabel $\beta$ [l] at 199 67
\endlabellist
\centering
\includegraphics[width=0.35\textwidth]{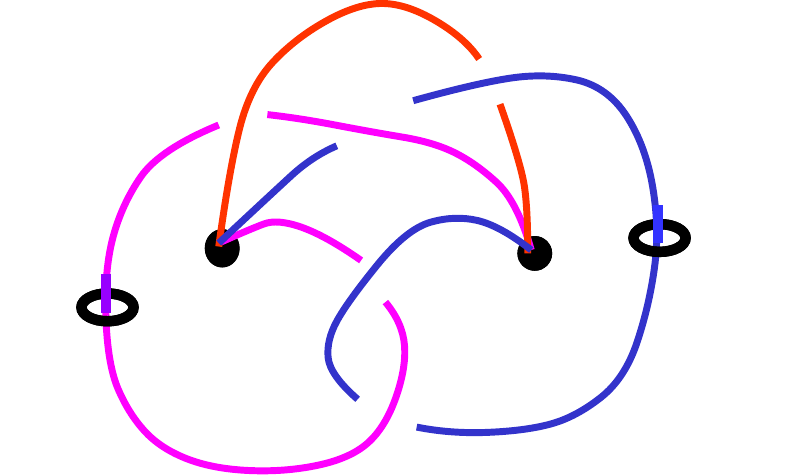}
\caption{} \label{fig:Kinoshita}
\end{figure}

\medskip
This naive conjecture is plainly false.  Most obviously, take $M$ to
be merely a regular neighborhood of $W$; no matter how $\aaa$ and
$\bbb$ are chosen, both $M[\aaa]$, $M[\bbb]$ are $\bdd$--reducible.
But there are easy counterexamples even for $M = S^3$.  For example,
attach an arc $e$ to the unknot $U \subset S^3$, an arc chosen to be
so complicated that the closed complement $S^3 - \eta(U \cup e)$ is
both irreducible and $\bdd$--irreducible.  Let $W = \eta(U \cup e)$
and choose both $\aaa$ and $\bbb$ to be copies of the meridian of $W$
that is dual to the arc $e$.  Then both $L[\aaa]$, $L[\bbb]$ are the
unknot.  Of course, taking $\aaa$ and $\bbb$ parallel like this might
be regarded as cheating.  \fullref{fig:Kinoshita} (due to Kinoshita
\cite{Ki}) is a more subtle counterexample in which $\aaa$ and $\bbb$
aren't parallel.  The complement of Kinoshita's graph is also called
the Thurston wye manifold and is known to be $\bdd$--irreducible, so
$W$ is knotted.

In view of these counterexamples to the naive conjecture, are there simple conditions that ensure the conclusion of the conjecture, that either $W$ is unknotted or at least one of $M[\aaa]$ or $M[\bbb]$ is both irreducible and $\bdd$--irreducible?  With reasonable conditions on the original pair $(M, W)$, it appears that there are such conditions on filling meridians $\aaa$ and $\bbb$.  These are outlined in the next section.

This research was partially supported by an NSF grant. Thanks also to Catalonia's Centre Recerca Matem\`atica for their extraordinary hospitality while this work was underway, and to both the referee and Scott Taylor for very helpful comments.

\section{A conjecture and a weak converse}

Suppose $W$ is a genus two handlebody properly embedded in a compact orientable $3$--manifold $M$.  The pair is called {\em admissible} if these conditions are satisfied:

\begin{itemize}
\item  Any sphere in $M$ is separating.
\item $M$ contains no Lens space summands.
\item Any pair of curves in $\bdd M$ that compress in $M$ are isotopic in $\bdd M$.
\item $M - W$ is irreducible.
\item $\bdd M$ is incompressible in the complement of $W$.
\end{itemize}

These are reasonable conditions to assume in our context:  The first two guarantee that the complement of a knot in a closed $M$ is $\bdd$--reducible only if the knot is the unknot.  The third condition (which is the most technical) removes the first  counterexample above, in which $M$ is merely a regular neighborhood of $W$.  The last two conditions remove obvious counterexamples in which reducing spheres or $\bdd$--reducing disks exist even before filling meridians are added.

The precise conditions that we propose on the pair of filling meridians $\aaa, \bbb \subset W$ depend on whether $\aaa$ and $\bbb$ are separating or non-separating.  Suppose that $\aaa$ and $\bbb$ have been properly isotoped in $W$ to minimize $|\aaa \cap \bbb|$.  In particular $\aaa \cap \bbb$ consists of a possibly empty collection of arcs.  

\begin{defin} \label{defin:aligned1} If $\aaa$ and $\bbb$ are both non-separating then they are {\it aligned\/} if 

\begin{itemize}

\item there is a non-separating meridian disk for $W$ that is disjoint from both $\aaa$ and $\bbb$ and

\item all arcs of $\aaa \cap \bbb$ are parallel in both disks.

\end{itemize}

\end{defin}  

\begin{defin} \label{defin:aligned2} If $\aaa$ and $\bbb$ are both separating then they are {\it aligned\/} if 

\begin{itemize}

\item there is a non-separating meridian disk $\mu$ for $W$ that is disjoint from both $\aaa$ and $\bbb$ and

\item there is a longitude in the boundary of the solid torus $W - \mu$ that is disjoint from $\aaa$ and $\bbb$.

\end{itemize}

\end{defin}

Here a {\em longitude} of a solid torus means any simple closed curve in the boundary that intersects a meridian disk in a single point.  If $\aaa$ and $\bbb$ are aligned and separating, then they both lie in the $4$--punctured sphere $\bdd W - (\bdd \mu \cup \mathrm{longitude})$ and they separate the same pairs of punctures there.  

\begin{figure}[ht!]
\labellist\small
\pinlabel $\mu$ [b] at 122 163
\pinlabel $\alpha$ at 203 18
\hair 1pt
\pinlabel $\beta$ [l] at 188 144
\pinlabel $\beta'$ [b] at 53 135
\endlabellist
\centering
\includegraphics[width=0.4\textwidth]{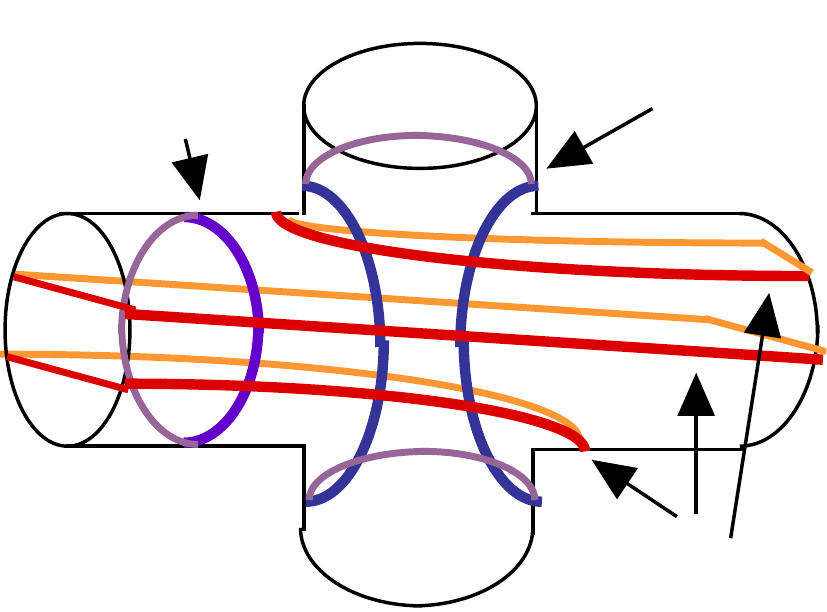}
\caption{} \label{fig:aligned}
\end{figure}

\begin{defin} \label{defin:aligned3}  If $\aaa$ is non-separating and $\bbb$ is separating, then the disks are {\it aligned\/} if they are disjoint or, when they are not disjoint,
\begin{itemize} 
\item one solid torus component of $W - \bbb$ has a meridian $\mu$
disjoint from both $\aaa$ and $\bbb$ and
\item the other solid torus component of $W - \bbb$ has a meridian $\bbb'$ that is disjoint from $\bbb$ and furthermore
\item $\bbb'$ is maximally aligned with $\aaa$.  That is, $\aaa$ and $\bbb'$ are aligned (as defined in \fullref{defin:aligned1}) and $|\bbb' \cap \aaa| = |\bbb \cap \aaa| - 1$ .  See \fullref{fig:aligned}.
\end{itemize}
\end{defin}  

For a pair of non-separating disks, the condition that all arcs of   $\aaa \cap \bbb$ are parallel in both disks means that either $\aaa$ and $\bbb$ are disjoint, or exactly two components of $\aaa - \bbb$ (and two components of $\bbb - \aaa$) are outermost disks, that is disks incident to a single arc component of $\aaa \cap \bbb$.  The condition is clearly satisfied whenever there are at most two arcs of intersection, ie if $|\bdd \aaa \cap \bdd \bbb| \leq 4$.

Here is the main conjecture:

\begin{conj} \label{conj:main}   
If $(M, W)$ is an admissible pair then

\begin{itemize}

\item $W$ is unknotted and $M = S^3$, or

\item at least one of $M[\aaa]$ or $M[\bbb]$  is both irreducible and $\bdd$--irreducible, or

\item $\aaa$ and $\bbb$ are aligned in $M$.

\end{itemize}    

\end{conj}

In many appearances of this general problem cited in the literature above, $\bdd \aaa$ and $\bdd \bbb$ intersect in few points, so they are automatically aligned.  In any case, the conjecture would only partially recapture the known results but would extend them to a more general setting.  \fullref{conj:main} appears to be true, at least with one additional technical hypothesis: When $|\bdd \aaa \cap \bdd \bbb| \leq 4$, assume further that $M$ contains no proper summand that is a rational homology sphere. A complete proof is not yet written down; even for a weaker result, in which $\bdd$--irreducible is removed from the second conclusion, the combinatorial argument is extremely complicated.  The intention here is to offer the more straightforward proofs in these three important special cases:
\begin{itemize}
\item  $M - W$ is $\bdd$--reducible.
\item $|\bdd \aaa \cap \bdd \bbb| \leq 4$\qua (This requires the additional technical condition.)
\item Both $\aaa$ and $\bbb$ are separating.
\end{itemize}
In addition, we explain why the combinatorics becomes so difficult once non-separating meridians are considered.

Before starting to verify the conjecture in these special cases, here is a sort of weak converse to \fullref{conj:main}. The assumption in \fullref{thm:converse}  that 
$\aaa$ and $\bbb$ are aligned is meant to be the weakest natural condition for which 
the theorem is likely to be true.

\begin{thm} \label{thm:converse}   
Suppose $\aaa$ and $\bbb$ are aligned meridians in $W$.  Then there is an \underline{unknotted} embedding of $W$ in $S^3$ (hence in any $3$--manifold) so that each of $L[\aaa]$ and $L[\bbb]$ is either the unknot or a split link.
\end{thm}

\begin{proof}  Suppose first that $\aaa$ and $\bbb$ are both separating.   Embed $W$ as the regular neighborhood of an eyeglass graph (ie two circles $\sss_1, \sss_2$  connected by an edge $e$) in $S^2 \subset S^3$. Since $\aaa$ is separating, such an embedding of $W$ can be found so that $\aaa$ is the meridian dual to the edge $e$.  Then  $L[\aaa] = \sss_1 \cup \sss_2$ is split.  Further choose the embedding $W$ so that in the framing of the tori $W - \aaa$, the longitude $\lambda$ of the solid torus $W - \mu$ that is disjoint from $\aaa \cup \bbb$ (cf \fullref{defin:aligned2}) is one of the three curves $\bdd W \cap S^2$.  Then $\lambda$ bounds a disk in $S^2 - W$, so one component of $L[\bbb]$ bounds a disk in the complement of the other, showing that the link $L[\bbb]$ is also  split.

In case $\aaa$ and $\bbb$ are non-separating, a more subtle construction is required.  Begin with an annulus $A$ in $S^2 \subset S^3$ and draw a pair $P$ of disjoint spanning arcs.   Let $\mu_{\pm}$ be points in $\bdd A - \bdd P$, one in each component of $\bdd A$ but both lying on the boundary of the same rectangle component of $A - P$.  A product neighborhood $A \times [-1, 1] \subset S^2 \times [-1, 1] \subset S^3$ is a solid unknotted torus in which $P \times [-1, 1] $ can be thought of as the union of two disjoint meridians $\aaa$ and $\bbb$.  Connect disk neighborhoods of $\mu_{\pm} \times \{ 1/2 \}$ in $\bdd A \times [-1, 1] $ by adding an unknotted $1$-handle on the outside of $A \times [-1, 1]$ that lies above $S^2 \times \{ 0 \}$.  The result is an unknotted embedding of $W$, with the meridian $\mu$ of the $1$--handle also a meridian of $W$ that is disjoint from $\aaa$ and $\bbb$.   Now repeat the construction after first altering exactly one of the original spanning arcs in $P$ by $n  \geq 0$ Dehn twists around the core of the annulus $A$.  The construction gives aligned meridians with exactly $n - 1$ arcs of intersection.  It is not hard to show that, given a pair of aligned non-separating meridians in $W$ which have $n - 1$ arcs of intersection, there is an automorphism of $W$ that carries the pair to $\aaa$ and $\bbb$.  The corresponding knots $L[\aaa]$ and $L[\bbb]$ are easily seen to be unknotted.   (In fact they are isotopic: the isotopy merely undoes the Dehn twists used in the construction by  adding twists around the meridian $\mu$ of the $1$--handle.)  See \fullref{fig:converse1}.

\begin{figure}[ht!]
\labellist\small
\pinlabel $\mu_+$ [b] at 235 385
\pinlabel $\mu_-$ [t] at 230 345
\pinlabel $\mu$ [t] <0pt,2pt> at 326 94
\pinlabel $A$ [r] at 106 314
\pinlabel $A\times\{-1\}$ [tr] at 66 109
\pinlabel $A\times\{1\}$ [t] at 113 76
\pinlabel $\alpha$ [b] at 164 129
\pinlabel $\alpha$ [b] at 164 96
\pinlabel $\alpha$ [b] at 413 136
\pinlabel* $\beta$ [br] at 171 205
\endlabellist
\centering
\includegraphics[width=0.7\textwidth]{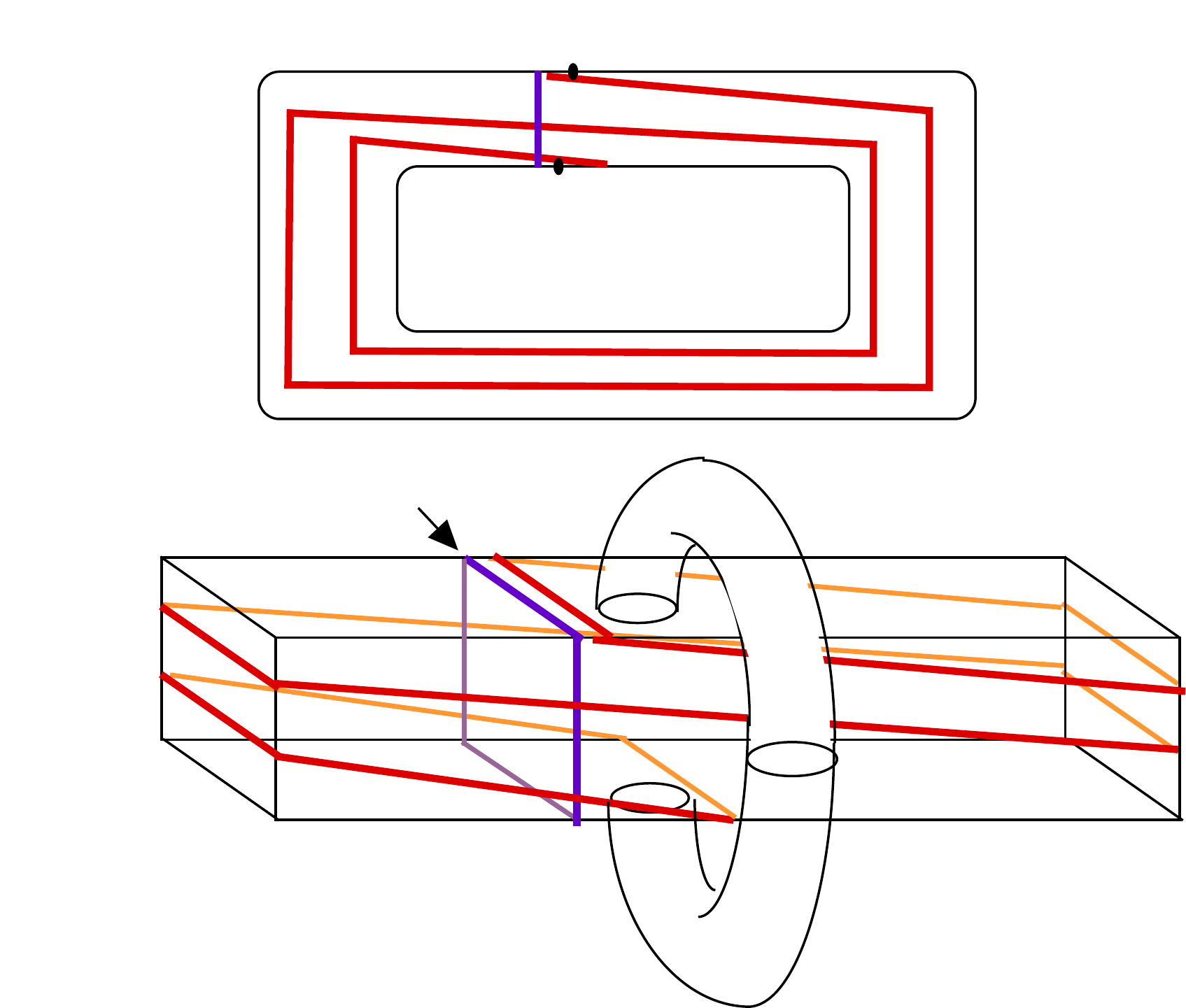}
\caption{} \label{fig:converse1}
\end{figure}

In case $\bbb$ is separating and $\aaa$ is not, let $\mu$ and $\bbb'$ be meridians of the solid tori $M - \bbb$ as given in \fullref{defin:aligned3}, so that 
\begin{itemize}
\item both $\mu$ and $\bbb'$ are disjoint from $\bbb$, 
\item $\mu$ is also disjoint from $\aaa$ and 
\item $\bbb'$ is aligned with $\aaa$.
\end{itemize}
Apply the previous construction to the aligned meridians $\aaa$ and $\bbb'$.  The definition requires that the aligned $\bbb'$ intersect $\aaa$ in almost as many components as $\bbb$ does.  Viewed dually, this implies that $\bbb$ intersects $\aaa$ only once more than $\bbb'$ does.  This constrains $\bbb$ to be a regular neighborhood of the arc $\bbb_-$ connecting  $\mu_+ \times \{ 1/2 \}$ to $\mu_-  \times \{ 1/2 \}$ in the disk $\bdd(A \times I) - (A \times \{ -1 \} \cup \bdd \bbb')$. As before, $W - \aaa$ is the unknot.   $W - \bbb$ is a regular neighborhood of the trivial link, with one component parallel to the core of $A \times \{ 0 \}$ and the other to the union of $\bbb_-$ and the core of the $1$--handle. See \fullref{fig:converse2}.
\end{proof}

\begin{figure}[ht!]
\labellist\small
\pinlabel* $\mu$ [t] at 261 376
{\hair 1.5pt \pinlabel $\mu$ [t] at 274 95 }
\pinlabel $L_\beta$ <0pt, 1pt> at 120 226
\pinlabel $\alpha$ [b] at 120 376
\pinlabel $\alpha$ [b] at 120 412
\pinlabel $\alpha$ [b] at 352 419
\pinlabel* $\beta$ [br] <2pt, 0pt> at 161 490
\endlabellist
\centering
\includegraphics[width=0.6\textwidth]{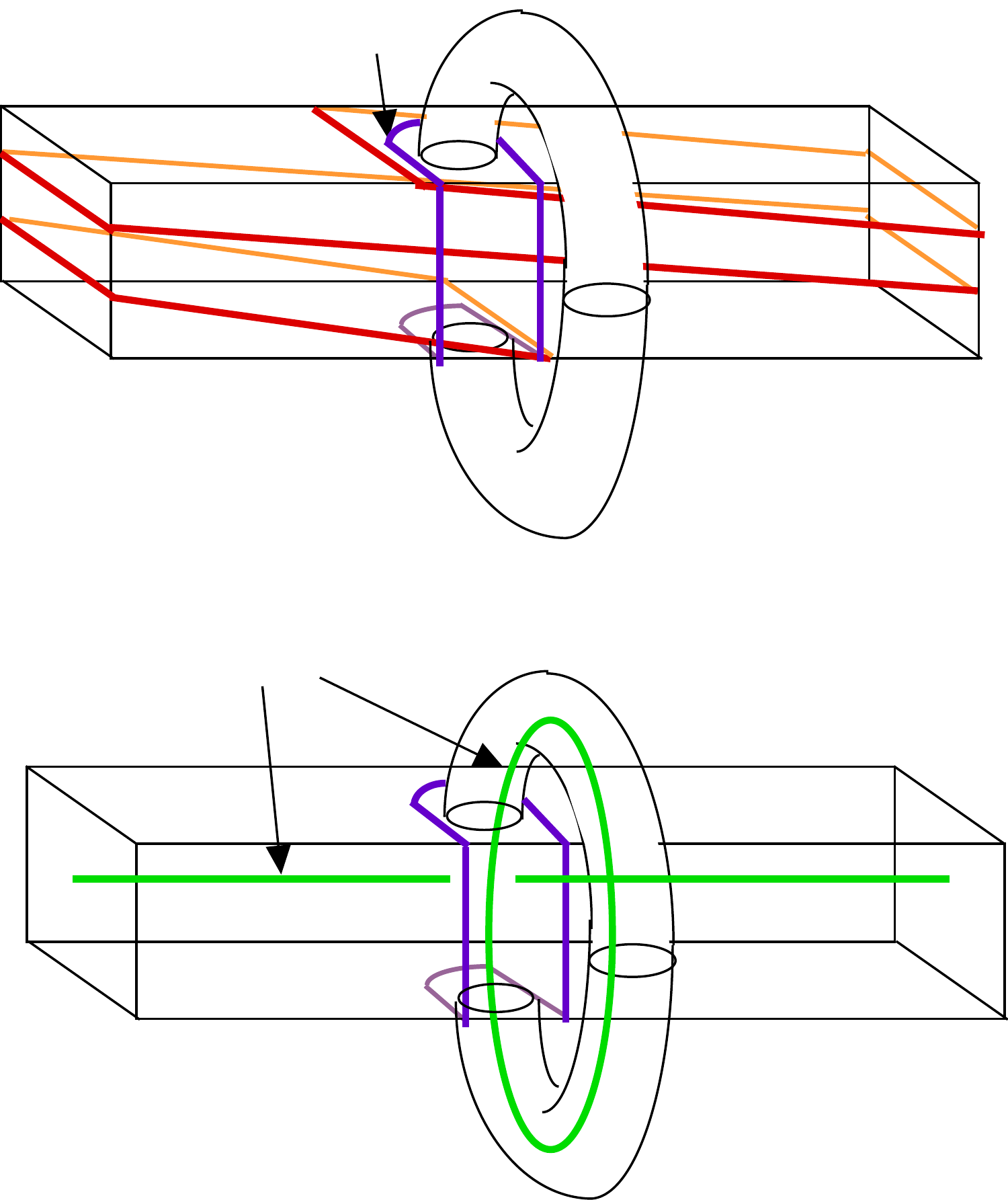}
\caption{} \label{fig:converse2}
\end{figure}  

A first step towards the proof of \fullref{conj:main} is to note that we can restrict to the case in which $M - W$ is $\bdd$--irreducible. 

\begin{prop} \label{prop:mainchange1}   \fullref{conj:main} is true when $M - W$ is $\bdd$--reducible.
\end{prop}

\begin{proof}\
{\bf Case 1}\qua There is a $\bdd$--reducing disk $P$ for $M - W$ such that $|\bdd P \cap \bdd \aaa| = 1$ (or, symmetrically, $|\bdd P \cap \bdd \bbb| = 1$).

In this case, $W - \eta(\aaa) \cong W \cup \eta(P)$ so any reducing sphere for $\Maa$ would be a reducing sphere for $M - W \supset M - (W \cup \eta(P))$, contradicting the assumption that $(M, W)$ is an admissible pair.  Similarly, if $\Maa$ is $\bdd$--reducible via a disk whose boundary lies on $\bdd M$, the disk would be a $\bdd$--reducing disk for $M - W$, contradicting the assumption that $(M, W)$ is an admissible pair.  Hence either $\Maa$ is both irreducible and $\bdd$--irreducible (the second conclusion) or there is a $\bdd$--reducing disk for $\Maa$ whose boundary lies on $\bdd (W - \eta(\aaa))$.  Since $M$ contains no Lens space summands, such a disk is incident to the solid torus $W - \eta(\aaa)$ in a longitude and the union of the disk with $W - \eta(\aaa)$ is a ball. Equivalently, there is a $\bdd$--reducing disk $D$ for $M - (W \cup \eta(P))$ so that $W \cup \eta(P) \cup \eta(D)$ is a ball.  $W$ is clearly a genus $2$ Heegaard surface for the ball, so  it follows from Waldhausen's theorem \cite{Wa} that $W$ is unknotted, giving the first conclusion and completing the proof in this case.  

\medskip
 {\bf Case 2}\qua There is a $\bdd$--reducing disk $(P, \bdd P) \subset (M - W, \bdd W - \aaa)$ for $M[\aaa]$ (or symmetrically for $\bbb$ and $\Mbb$).

\medskip
In this case first note that $W - \eta(\aaa)$ consists of one or two solid tori (depending on whether $\aaa$ separates $W$) and $\bdd P$ lies on one of them.  Since $M$ contains no Lens space or non-separating $2$--spheres, $\bdd P$ is in fact a longitude of the solid torus on whose boundary it lies.  So $\bdd P$ is non-separating.  

Suppose first that $\bdd P$ is disjoint from $\bdd \bbb$ as well as $\bdd \aaa$; we show that $\aaa$ and $\bbb$ are then aligned.  Let $\mu$ be a meridian for $W$ which is disjoint from $\aaa$ and whose boundary is disjoint from $\bdd P$.  In particular, if $\aaa$ is non-separating, take a parallel copy of $\aaa$ for $\mu$. If $\bbb$ intersects $\mu$ then an outermost disk of $\bbb$ cut off by $\mu$ would be a meridian of the solid torus $W - \mu$ for which $\bdd P$ is a longitude, a contradiction.  Hence $\bbb$ lies entirely in the solid torus $W - \mu$.  $\bbb$ can't be essential in that solid torus for the same reason.  Hence $\bbb$ is inessential in $W - \mu$ and so is either parallel to $\mu$ (hence disjoint from and therefore aligned with $\aaa$) or separating. If $\aaa$ is non-separating then it is parallel to $\mu$ and so disjoint from and therefore aligned with $\bbb$.  If $\aaa$, like $\bbb$, is separating, then $\bdd P$ is the longitude required by \fullref{defin:aligned2} to show that $\aaa$ and $\bbb$ are aligned.  

So henceforth assume that $\bdd P$ and $\bdd \bbb$ are not disjoint.
Apply the Jaco handle-addition lemma \cite{Ja} (as generalized by
Casson and Gordon \cite{CG} to the reducible case) to the $2$--handle
$\eta(\bbb)$ attached to the $\bdd$--reducible (via $P$) manifold
$M-W$: If $\Mbb$ is reducible or $\bdd$--reducible then there is a
$\bdd$--reducing disk $J$ for $M-W$ whose boundary is disjoint from
$\bbb$.  Since $\bdd M$ is incompressible in $M-W$, $\bdd J$ lies on
$\bdd W$.  Assume that such a disk $J$ has been chosen to minimize $|J
\cap P|$.

Suppose first that $J \cap P \neq \emptyset$ and consider an outermost disk $E$ cut off from $J$ by $P$.  We can assume the arc $\bdd E - \bdd P$ is essential in the surface $\bdd W - \bdd P$, else an isotopy of that arc across $\bdd P$ would reduce $|J \cap P|$.  It follows that $\bdd E$ is essential on the torus $\bdd(W \cup \eta(P))$, for the alternative is that it cuts off a disk containing all of one copy of $P$ in the torus but only part of the other, and so, absurdly, $\bdd J$ intersects one side of $\bdd P$ more often than the other.  Then $W \cup \eta(P) \cup \eta(E)$ is the union of a solid torus and an essential $2$--handle.  Since $M$ contains neither non-separating spheres nor Lens spaces,  $W \cup \eta(P) \cup \eta(E)$ is a $3$--ball.  Cap the $3$--ball off with a $3$--handle and the result is a genus two Heegaard splitting of $S^3$, in which $W$ is one of the two handlebodies.  But Waldhausen's theorem \cite{Wa} states that any Heegaard splittings of $S^3$ is isotopic to a standard one.  It follows that $W$ may be isotoped in the $3$--sphere (hence in a neighborhood of the $3$--ball $W \cup \eta(P) \cup \eta(E)$) so that it is standard.  Thus $W$ is unknotted in $M$ and we are done.  

The remaining case is that $J$ is disjoint from $P$.  If $\bdd J$ is essential on the torus  $\bdd(W \cup \eta(P))$ then use $J$ for $E$ in the previous argument and again we are done.   Suppose then that $\bdd J$ is inessential on $\bdd(W \cup \eta(P))$.  The disk it bounds can't contain just one copy of $P$, since $\bdd P$ intersects $\bbb$ and $\bdd J$ does not.   So the disk contains both copies of $P$, hence $\bdd J$ is separating in $\bdd W$.  Since $\bdd J$ is disjoint from the meridian $\bbb$,  $\bdd J$ then also bounds a separating meridian $J' \subset W$.  Since $\bdd \bbb$ and $\bdd P$ intersect, $\bbb$ is a meridian disk for the solid torus component of $W - J'$ for which $\bdd P$ is the longitude.  Hence $|\bdd \bbb \cap \bdd P| = 1$ and we are done by the first case.

Now consider the general case:

Apply the Jaco handle-addition lemma to $\eta(\aaa)$ and conclude that either $\Maa$ is irreducible and $\bdd$--irreducible (the second conclusion) or a $\bdd$--reducing disk $(P, \bdd P)$ can be found whose boundary is disjoint from $\aaa$.  Since $\bdd M$ is incompressible in $M - W$, $\bdd P \subset \bdd W - \aaa$.  If $\bdd P$ is essential on $\bdd \Maa$ the proof follows from Case 2.  If $\bdd P$ is inessential on $\bdd \Maa$ then it is coplanar (possibly parallel with) $\aaa$, so $\bdd P$ bounds a meridian $P'$ for $W$.  In particular, $\bdd P$ is separating, since $M$ contains no non-separating spheres.  Similarly there is a separating disk $(Q, \bdd Q) \subset (M - W, \bdd W - \bbb)$.  If $\bdd P$ and $\bdd Q$ are parallel in $\bdd W$ then it follows easily that $\aaa$ and $\bbb$ are disjoint from each other and hence aligned.  If $\bdd P$ and $\bdd Q$ are not isotopic in $\bdd W$ then they must intersect, and an outermost arc of intersection in $Q$ cuts off a disk $(F, \bdd F)  \subset (M-W, \bdd W - \bdd P)$ whose boundary is essential on, and hence a longitude of, one of the solid tori $W - P'$.  If $\aaa$ is a meridian of that solid torus, use $F$ as $P$ in Case 1.  If $\aaa$ is the meridian of the other solid torus, use $F$ as $P$ in Case 2.
\end{proof}

Following \fullref{prop:mainchange1}, \fullref{conj:main} is equivalent to the apparently weaker conjecture:

\begin{conj} \label{conj:main2}   
Suppose $(M, W)$ is an admissible pair and $M - W$ is $\bdd$--irreducible.
Then either 
\begin{itemize}
\item at least one of $M[\aaa]$ or $M[\bbb]$  is both irreducible and $\bdd$--irreducible, or
\item $\aaa$ and $\bbb$ are aligned in $M$.
\end{itemize}    
\end{conj}

\begin{prop} \label{prop:leq4}

When $|\bdd \aaa \cap \bdd \bbb | \leq 4$, \fullref{conj:main2} is true .  If, furthermore, $\aaa$ is separating and $M$ contains no proper summand that is a rational homology sphere then one of $\Maa$ or $\Mbb$ is irreducible and $\bdd$--irreducible.

\end{prop}

\begin{proof} When $|\bdd \aaa \cap \bdd \bbb | \leq 4$, there are at most two arcs of intersection, so each arc of intersection is incident to an outermost disk in both $\aaa$ and $\bbb$ and, moreover, all arcs of intersection are parallel in both disks.  Gluing an outermost disk of one to an outermost disk of the other gives a meridian disk for $W$ that is disjoint from both $\aaa$ and $\bbb$.  It follows that if $\aaa$ and $\bbb$ are both non-separating, then they are aligned.  

So it suffices to consider the case in which at least one of them, say $\aaa$, is separating.  We aim to show that the disks are aligned and, assuming $M$ contains no proper summand that is a rational homology sphere,  if $\Maa$ is either reducible or $\bdd$--reducible, then $\Mbb$ is neither reducible nor $\bdd$--reducible. There are two cases to consider:

\medskip

{\bf Case 1}\qua $|\bdd \aaa \cap \bdd \bbb | = 2$ 

\medskip

In this case $\bbb$ must be non-separating and it is immediate from
\fullref{defin:aligned3} that $\aaa$ and $\bbb$ are aligned.
The argument that one of $\Maa$ or $\Mbb$ is irreducible and
$\bdd$--irreducible now mimics in this more general setting the
central argument of Gabai \cite{Ga}, which in fact provides the
complete argument when $M$ is the $3$--sphere.

Let $L \subset M - W$ be the regular neighborhood of a circle parallel to $\bdd \aaa$, so a longitude of $\bdd L$ bounds a disk $D$ in $M - L$ that intersects $W$ in a single copy of $\aaa$.  Since $\aaa \cap \bbb$ is a single arc, $D$ intersects the solid torus $T = W - \bbb$ in two oppositely-oriented meridians.  (For example, in Gabai's setting, $W - \aaa$ is a split link and $W - \bbb$ is obtained from the split link by a band sum.  $\bdd D$ encircles the band.)  

First observe that $\Mbb - L$ is irreducible and $\bdd$--irreducible.  For suppose first that a $\bdd$--reducing disk or reducing sphere $\hat{Q}$ were disjoint from $D$.  Since $\Mbb \cup \eta(D) \cong W$,  $\hat{Q}$ would also be a $\bdd$--reducing disk or reducing sphere in $M - W$, contradicting the hypothesis.  Similarly, if $\hat{Q}$ intersects $D$ an innermost disk $E$ in $\hat{Q} - D$ provides a $\bdd$--reducing disk for $\bdd W$ in the complement of $\Mbb \cup \eta(D) \cong W$. 

Let $\Mbb_0$ be the manifold obtained from $\Mbb$ by $0$--framed surgery on $L$, that is, surgery with slope $\bdd D$.  

\medskip
{\bf Claim}\qua  {\sl If $\Maa$ is reducible or $\bdd$--reducible, then so is $\Mbb_0$.}

\medskip
To prove the claim, suppose $\hat{P}$ is a reducing sphere or $\bdd$--reducing disk for $\Maa$.  Then $\inter(\hat{P}) \cap W$ is a collection of parallel copies of $\aaa$.  Let $A \subset \bdd W$ be an annulus containing all their boundary components $\inter(\hat{P}) \cap \bdd W$.  Let $A'$ be a copy of $A$ pushed into $W$ rel $\bdd A$.  Then $W - \eta(A')$ is isotopic to $W \cup L$ and under this isotopy the boundary components of the planar surface $P = \hat{P}  - W$ that previously were on $A$ become $0$--framed curves on $\bdd (L)$.  After $0$--framed surgery on $L$ to get $\Mbb_0$, these components on $\bdd L$  can be capped off to give a copy $\hat{P}'  \subset \Mbb_0$ of $\hat{P}$ that either $\bdd$--reduces or reduces $\Mbb_0$. See \fullref{fig:leq4}.  That $\hat{P}'$ does not bound a ball is obvious if $\hat{P}'$ is non-separating (eg $m$ is odd).  If $\hat{P}'$ is separating, note that such a ball must have come from the component of $M - (W \cup L \cup \eta(P))$ not adjacent to $W$.  But this component is completed by attaching $2$--handles in the same way in both $\Mbb_0$ and $\Maa$, and we know that $\hat{P}$ is a reducing sphere in $\Maa$ and so does not bound a ball. This completes the proof of the Claim.

\begin{figure}[ht!]
\labellist\small
\pinlabel $\eta(L)$ [t] at 118 25
\pinlabel $A$  at 141 231
\pinlabel $A'$ [t] at 194 30
\pinlabel $\alpha_1\ldots\alpha_m$ [b] at 139 340
\pinlabel $\alpha_1\ldots\alpha_m$ [b] at 148 130
\pinlabel $\beta$ [l] <0pt,2pt> at 231 302
\endlabellist
\centering
\includegraphics[width=0.45\textwidth]{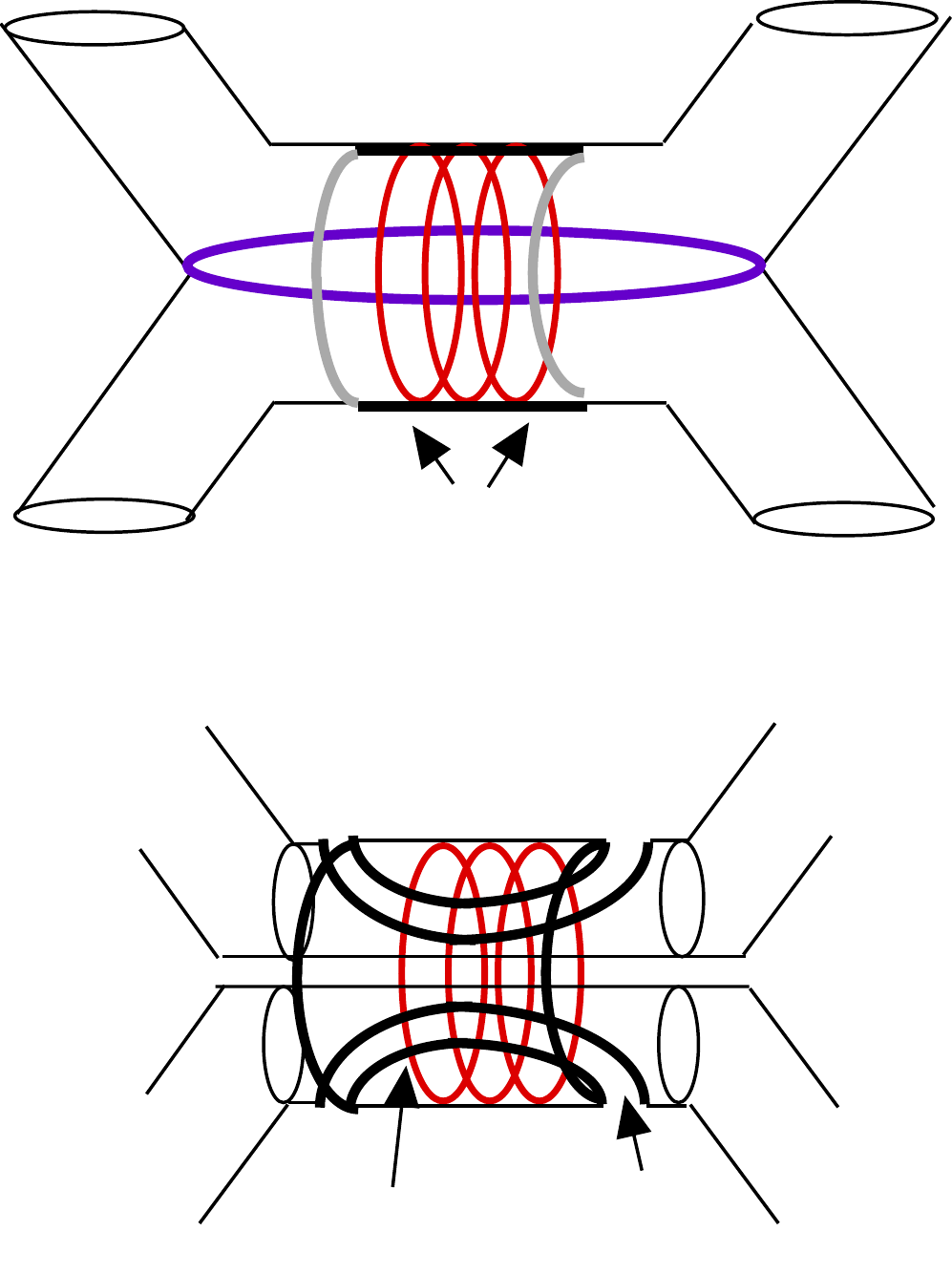}
\caption{} \label{fig:leq4}
\end{figure}  

\medskip
Following the Claim, apply the central theorem of Scharlemann
\cite{Sc4} to the link $L$ in the manifold $\Mbb$ and rule out
possible conclusions a)--c): $L \subset \Mbb$ is not homeomorphic to a
braid in $D^2 \times S^1$ since, for one thing, $L$ is
null-homologous.  Nor does $\Mbb$ contain a Lens space summand, since
$M$ doesn't.  A longitude of $\bdd L$ bounds a disk in $M - L$ so
$\bdd L$ can't be cabled in $M$ with the same boundary slope.  Hence
conclusion d) of the main theorem of \cite{Sc4} applies to $\Mbb$ and
$\Mbb_0$.  We now examine the consequences.

Suppose that $\Maa$ is reducible.  Then $\Mbb_0$ is reducible, so by  \cite{Sc4}, $\Mbb$ is irreducible and $\bdd$--irreducible.  

Suppose that $\Maa$ is $\bdd$--reducible via the disk $\hat{P}$.  We note that if $\bdd \hat{P}$ lies on the boundary of one of the solid tori components of $W - \aaa$ then that solid torus lies in a $3$--ball in $M$, so $\Maa$ is reducible.  Then, via the previous comment, $\Mbb$ is irreducible and $\bdd$--irreducible. So we can assume that $\bdd \hat{P} \subset \bdd M$.  Then $\bdd \hat{P}$ is disjoint from any $\bdd$--reducing disk for $\Mbb$ since, by hypothesis, any essential simple closed curves in $\bdd M$ that compress in $M$ are isotopic.  Hence $\Mbb_0$ is $\bdd$--reducible via a disk whose boundary is disjoint from any $\bdd$--reducing disk for $\Mbb$.  So by  \cite{Sc4} $\Mbb$ is irreducible and $\bdd$--irreducible.   To summarize: If $\Maa$ is reducible or $\bdd$--reducible, then $\Mbb$ is neither reducible nor $\bdd$--reducible.  Thus we have the first conclusion of \fullref{conj:main2},  completing the proof in this case.  

\medskip
{\bf Case 2}\qua  $|\bdd \aaa \cap \bdd \bbb | = 4$ 

\medskip
In this case $\bbb \cap \aaa$ consists of two arcs; call the rectangles that lie between them in $\aaa$ and $\bbb$ respectively $R_{\aaa}$ and $R_{\bbb}$.  The outermost disks of $\bbb$ cut off by $\aaa$ are both meridians of the same solid torus component $W_1$ of $W - \aaa$; in fact they are parallel copies of the same disk meridian $D_{\bbb} \subset W_1$ whose boundary consists of an arc in $R_{\aaa}$ parallel to the arcs $\bbb \cap \aaa$ and an arc on $\bdd W - \aaa$.  See \fullref{fig:leq42}.  The meridian $D_{\bbb}$ is properly isotopic in $W$ to a meridian that is disjoint from both $\aaa$ and $\bbb$; a meridian for the other solid torus component of $W - \aaa$ intersects $\bbb$ in a single arc.  It follows immediately that $\aaa$ and $\bbb$ are aligned.   It remains to show that if $M$ contains no proper summand that is a rational homology sphere and if $\Maa$ is reducible or $\bdd$--reducible, then $\Mbb$ is neither reducible nor $\bdd$--reducible.

\begin{figure}[ht!]
\labellist\small\hair 1pt
\pinlabel $W'$ [tl] at 291 362
\pinlabel $W_2$ [t] at 150 47
\pinlabel $D_\beta$ [b] at 136 153
\pinlabel $\alpha$ [t] at 149 309
\pinlabel* $\alpha'$ [t] <1pt,0pt> at 165 114
\pinlabel $\beta$ [l] at 282 261
\endlabellist
\centering
\includegraphics[width=0.45\textwidth]{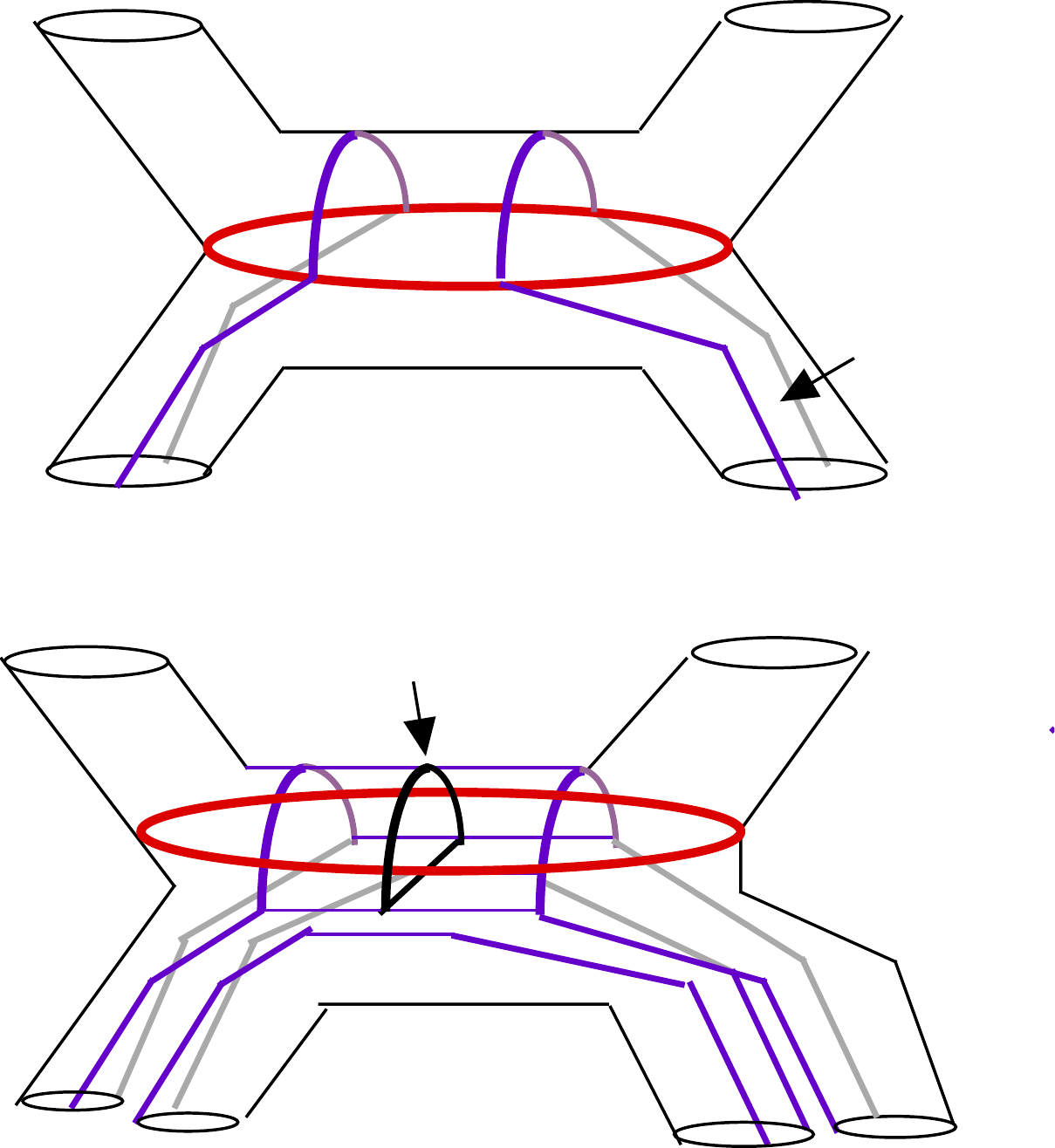}
\caption{} \label{fig:leq42}
\end{figure}

 The union $A = R_{\aaa} \cup R_{\bbb}$ is either an annulus or M\"obius band in $W$ (depending on the orientations of the arcs $\bbb \cap \aaa$ in $\aaa$) which, when $\bdd$--compressed to $\bdd W$ along $D_{\bbb}$ becomes $\bbb$.  Put another way, $W_- = W - A$ is a handlebody or pair of handlebodies (depending on whether $A$ is a M\"obius band or an annulus) which, when $\bdd$--reduced along $D_{\bbb}$ becomes $W - \bbb$.   A copy $\aaa'$ of $\aaa$ pushed into $W_1$ intersects $D_{\bbb}$ in a single arc.  Essentially, we intend to settle Case 2 by applying Case 1 to the disks $D_{\bbb}$ and $\aaa'$, two meridians of $W_-$ that intersect in a single arc. 

If $A$ is an annulus, then $W_-$ is the union of a solid torus $W_2$ (essentially the component of $W - \aaa$ that is not $W_1$) and a genus two component.  In this case, let $W'$ be the genus $2$--component and $M' = M - W_2$.   If $A$ is a M\"obius band let $W' = W_-$ and $M' = M$.  

\medskip
{\bf Claim 1}\qua  {\sl If $\Maa$ is reducible or $\bdd$--reducible, 
so is $M'[\aaa]$.} 

\medskip
{\bf Proof of Claim 1}\qua $M'[\aaa]$ is obtained from $\Maa$ by attaching $A$ to $\bdd \Maa$ via $\bdd A$ and thickening $A$.  This operation can't turn a complementary component of a sphere or properly embedded disk into a $3$--ball, so it can't make a reducible manifold irreducible or a $\bdd$--reducible manifold $\bdd$--irreducible.  

\medskip
{\bf Claim 2}\qua  {\sl If $\Mbb$ is reducible or $\bdd$--reducible, so is $M'[D_{\bbb}]$.}  

\medskip
{\bf Proof of Claim 2}\qua  By construction, $\bdd$--compressing $A$ to $\bdd W$ via $D_{\bbb}$ gives $\bbb$, so $\Mbb$ is homeomorphic to $M'[D_{\bbb}]$

\medskip
{\bf Claim 3}\qua  {\sl $M' - W'$ is irreducible and $\bdd$--irreducible.}  

\medskip
{\bf Proof of Claim 3}\qua Suppose $S$ is a reducing sphere or $\bdd$--reducing disk for $M' - W'$ and $A'$ is the annulus or M\"obius band dual to $A$ in $W$.  That is, $W' \cup \eta(A') \cong W$ if $A$ is a M\"obius band and $W' \cup \eta(A') \cup W_2 \cong W$ if $A$ is an annulus. If $S$ is disjoint from $A'$ then $S$ is a reducing sphere or $\bdd$--reducing disk for $M - W$, contradicting hypothesis.  

On the other hand, suppose $S$ and $A'$ are not disjoint.  Closed components of $S \cap A'$ that are inessential in $A'$ can be removed by a standard innermost disk argument; if any essential closed components of $S \cap A'$ remain, then an innermost disk of $S - A'$ would be a compressing disk for $\bdd W$ in $M- W$, again contradicting hypothesis.  Similarly, all arc components of $S \cap A'$ (which only can arise if $S$ is a $\bdd$--reducing disk, not a reducing sphere) that are inessential in $A'$ (ie non-spanning) can be removed by a standard outermost arc argument.  On the other hand, a spanning arc in $S \cap A'$ is clearly impossible if $A'$ is an annulus, since the ends of $A'$ are on different components of $\bdd(M' - W')$ whereas $\bdd S$ can only be incident to one.  If $A'$ is a M\"obius band then an outermost disk cut off from $S$ by $A'$ will be a $\bdd$--reducing disk for $M-W$, again contradicting hypothesis.

\medskip
Following Claim 3, Case 1 may be applied to the disks $\aaa, D_{\bbb} \subset W'$ in $M'$: If $M'[\aaa]$ is reducible or $\bdd$--reducible then $M'[D_{\bbb}]$ is irreducible and $\bdd$--irreducible.  Hence, following Claims 1 and 2, if $\Maa$ is reducible or $\bdd$--reducible, $\Mbb$ is irreducible and $\bdd$--irreducible. 
\end{proof}

\begin{prop} \label{prop:bothsep}   

If both $\aaa$ and $\bbb$ are separating \fullref{conj:main2} is true.  Furthermore, if, in addition to the assumptions of \fullref{conj:main2}, $|\bdd \aaa \cap \bdd \bbb | \geq 6$, then one of $\Maa$ or $\Mbb$ is irreducible and $\bdd$--irreducible.

\end{prop}

The rest of the section is devoted to the proof of \fullref{prop:bothsep}.  The proof requires a few internal lemmas.

Following \fullref{prop:leq4} we may as well assume $|\bdd \aaa \cap \bdd \bbb | \geq 6$. 

Since both $\aaa$ and $\bbb$ are separating, $W - \eta(\aaa)$ and $W - \eta(\bbb)$ are each a pair of solid tori.   If $\Maa$, say, is $\bdd$--reducible and the $\bdd$--reducing disk $P$ is incident to a solid torus component $W_1$ of $\bdd(W - \eta(\aaa))$ then $\bdd(\eta(W_1 \cup P))$ is a sphere in $M$ that separates the two solid tori, so $\Maa$ is reducible as well as $\bdd$--reducible.  On the other hand, if $\Maa$ is $\bdd$--reducible and the $\bdd$--reducing disk $P$ is incident to $\bdd M$ then $P$ is a $\bdd$--reducing disk for $M$ as well.  In that case $\bdd P$ lies in the unique isotopy class of $\bdd M$ which compresses in $M$, following the assumption that $(M, W)$ is admissible.  We deduce that if neither $\Maa$ nor $\Mbb$ are reducible but both are $\bdd$--reducible then the $\bdd$--reducing disks have disjoint boundaries, lying in $\bdd M$.  

The proof proceeds by assuming that both $\Maa$ and $\Mbb$ are $\bdd$--reducible or reducible and derives a contradiction by examining components of intersection between a reducing sphere or $\bdd$--reducing disk in $\Maa$ and a  reducing sphere or $\bdd$--reducing disk in $\Mbb$.  Following the comments above it would be natural to view the case of reducing spheres as entirely representative, since when there are no reducing spheres the $\bdd$--reducing disks that can be used necessarily have disjoint boundaries, so the combinatorics is unaffected by the boundary.  In fact we will ignore $\bdd$--reducing disks whose boundaries lie on $\bdd M$, because the combinatorics involved is identical to that for reducing spheres.  But for somewhat technical reasons, it is best to allow another possibility, even when there are reducing spheres.  Say that a $\bdd$--reducing disk for $\Maa$ is {\em special} if its boundary $p_0$ lies in $\bdd W$ and $p_0$ is parallel in $\bdd W$ to the end-point union of a subarc of $\bdd \aaa$ and a component $\bbb_0$ of $\bdd \bbb - \bdd \aaa$.  The curve $p_0$ is necessarily a longitude of the solid torus component $W_1$ of $W - \aaa$ on which it lies.  Furthermore, since each arc of $\bdd \bbb - \aaa$ in the punctured torus $\bdd W_1 - \aaa$ is essential and lies in the annulus $\bdd W_1 - (\aaa \cup \bbb_0)$, it follows that each such arc component intersects $p_0$ at most once.  There is a symmetric definition for a special $\bdd$--reducing disk for $\Mbb$. 

To begin the combinatorial argument, let $\hat{P}$ (resp $\hat{Q}$) be either a reducing sphere or a special $\bdd$--reducing disk for $\Maa$ (resp $\Mbb$).   Since $M-W$ is irreducible and $\bdd$--irreducible, $\hat{P}$ must pass through the handle $\eta(\aaa)$ $m \geq 1$ times. Among all reducing spheres and special $\bdd$--reducing disks for $\Maa$, choose $\hat{P}$ to minimize $m$.  Let $P = \hat{P} - W $ be the associated properly embedded planar surface in $M-W$.  Then $\bdd P$ has  components $\aaa_1,\ldots, \aaa_m, m \geq 1$, each of them parallel on $\bdd W$ to $\bdd \aaa$ and, if $\hat{P}$ was a special $\bdd$--reducing disk, a component $p_0 = \bdd \hat{P}$.  Label the $\aaa_i$ in the order they appear in an annular neighborhood of $\bdd \aaa$ in $\bdd W$. There are two choices for the direction of the ordering.  If $\hat{P}$ is a $\bdd$--reducing disk, choose the ordering which puts $\aaa_1$ adjacent to $p_0$; otherwise choose either ordering.  For $\hat{Q}$ there is a similar planar surface $(Q, \bdd Q) \subset (M - W, \bdd W)$ whose boundary components are similarly labeled $\bbb_1,\ldots\bbb_n, n \geq 1$ and possibly $q_0$.  Isotope $P$ and $Q$ so as to minimize $|P \cap Q|$.  $| P \cap Q| > 0$ since $\aaa$ and $\bbb$ intersect.   Let $\hat{\aaa}_i, \hat{\bbb}_j \subset W$ be meridian disks in $W$ bounded by $\aaa_i, \bbb_j$ respectively, for each $1 \leq i \leq m$ and $ 1 \leq j \leq n$.

As is now a classical strategy,  view $P \cap Q$ as giving rise to graphs $\Sss$ and $\Uuu$ respectively in $\hat{P}$ and $\hat{Q}$.  The vertices of the graphs correspond respectively to the disks $\hat{\aaa_i} \subset \hat{P}$ and $\hat{\bbb_i} \subset \hat{Q}$.  Edges of each graph correspond to the arcs of intersection in $P \cap Q$.  Circles of intersection are ignored.  The valence of each vertex in $\Sss$ is $|\bdd \aaa \cap \bdd \bbb| \cdot n + |\bdd \aaa \cap q_0|$ and in $\Uuu$ is $|\bdd \aaa \cap \bdd \bbb| \cdot m + |\bdd \bbb \cap p_0|$.  When $\hat{P}$ is a special $\bdd$--reducing disk, the number of ends of edges that are incident to $p_0 = \bdd \hat{P}$ is less than half the valence of each vertex of $\Sss$, since each arc of $\bdd Q - \aaa$ lying on the torus component of $\bdd (W - \aaa)$ containing $p_0$ contributes $2$ to the valence of $\aaa_1$ but intersects $p_0$ at most once, and at least one such component is disjoint from $p_0$.  Similarly, the number of ends of edges that are incident to $q_0 = \bdd \hat{Q}$ is less than half the valence of each vertex of $\Uuu$.

Consider a point $x$ in $\bdd P \cap \bdd Q$, say $x \in \aaa_i \cap \bbb_j$.   To $x$ assign the ordered pair $(i, j), 1 \leq i \leq m, 1 \leq j \leq n$.  When viewed in the graphs $\Sss$ and $\Uuu$ the point $x$ appears as an end of a unique edge.  Assign to the end of the edge in $\Sss$ the label $j$.  Similarly assign to the end of the edge in $\Uuu$ the label $i$.  For each $1 \leq j \leq n$,  of the ends of edges incident  to any vertex in $\Sss$, exactly $|\aaa \cap \bbb|$ will have label $j$.  A similar remark holds for labeling around a vertex of $\Uuu$.  Points (if any) in $p_0 \cap \bbb_j$, $\aaa_i \cap q_0$ and $p_0 \cap q_0$ similarly are assigned the ordered pairs $(0, j), (i, 0)$ and $(0,0)$ respectively, giving rise to ends of edges labeled $0$ in respectively $\Uuu$, $\Sss$, and both $\Sss$ and $\Uuu$.  

An important difference between the topology exploited in this proof
and that used in the analysis of Dehn surgery (cf Culler, Gordon,
Luecke and Shalen \cite{CGLS}) is that in the latter, any two
components $\aaa_i, \bbb_j$ of the boundaries of $P$ and $Q$
respectively always intersect in the torus boundary of the
$3$--manifold with the same orientation.  In the present situation,
each $\aaa_i$ will intersect each $\bbb_j$ with both orientations.
Indeed, since both bound meridian disks of $W$, the sum of the
orientations of all points of intersection between the two will be
trivial.

Let $A \subset \bdd W$ be the annulus between $\aaa_1$ and $\aaa_m$ whose core is $\bdd \aaa$.  To any proper arc in $A$ that spans $A$ and is properly isotoped to minimally intersect the $\aaa_i$, assign the orientation that points from $\aaa_1$ to $\aaa_m$.  This is called the {\em spanning} orientation of the arc.  Similarly let $B \subset \bdd W$ be the annulus between $\bbb_1$ and $\bbb_n$ and define the spanning orientation for spanning arcs of $B$ to be the orientation that points from $\bbb_1$ to $\bbb_n$.  The spanning orientation for $A$ gives what is called the {\em spanning} normal orientation in $\bdd W$ to each component $\aaa_1, \ldots, \aaa_m$ of $\bdd P$, namely the normal orientation that points from $\aaa_i$ to $\aaa_{i+1}$.  Similarly define the spanning normal orientation for each $\bbb_j \subset \bdd Q$.  

Fix normal orientations for $P$ and $Q$.  Following such a choice, an edge in $\Sss$ whose corresponding edge in $\Uuu$ connects two vertices (so neither end in $\Sss$ is labelled $0$) is called {\em incoherent} if the normal orientation of $Q$ along the arc of intersection agrees with the spanning orientation at one end and disagrees at the other.  Otherwise (if it agrees at both ends or disagrees at both ends) it is called {\em coherent}.  Whether the edge is coherent or incoherent is independent of the original choice of orderings of the $\aaa_i$ and the $\bbb_j$ or the choice of normal orientation for $Q$.  Note that by definition an incoherent edge cannot have ends with the same label $j$.  If two coherent edges are parallel and adjacent, one with ends labeled $j_1, j_1'$ and the other with corresponding ends labeled $j_2, j_2'$ then $j_1 - j_1' \equiv  j_2 - j_2'$ mod $n$.  Similarly if two incoherent edges are parallel and adjacent, then $j_1 + j_1' \equiv  j_2 + j_2'$ mod $n$.   Similar remarks hold for labels of ends of edges in $\Uuu$.  Note that if two edges are parallel and if one is coherent and the other is incoherent, then the spanning orientations agree on one end (ie, induce the same orientation on that component of $\bdd P$) and disagree on the other.  The only way that spanning orientations can disagree at adjacent labels is if the interval between them does not lie on $B$, so one pair of adjacent ends are either both labeled $1$ or both labeled~$n$. 

\begin{lemma}  \label{lemma:noloops} There are no trivial loops in either $\Uuu$ or $\Sss$.
\end{lemma}

\begin{proof}  We show that an innermost trivial loop in $\Uuu$ can be used to reduce $m$.  $W - \hat{P}$ consists of $m-1$ copies of $D^2 \times I$ (labelled $W_1,\ldots,W_{m-1}$) and two solid tori $W_0$ and $W_m$.  Each $1$--handle $W_i, 1 \leq i \leq m-1$ lies between copies $\hat{\aaa}_i$ and $\hat{\aaa}_{i+1}$ of $\aaa$; the solid torus $W_0$ is incident to $\hat{P}$ in $\hat{\aaa}_1$ (and possibly $p_0$) and the solid torus $W_m$ is incident to $\hat{P}$ in $\hat{\aaa}_m$.

An innermost trivial loop in $\Uuu$ cuts off a disk $D$ from $Q$ whose boundary consists of an arc $q$ on $\bdd Q - P \subset \bdd W - P$ (in the component of $\bdd Q$ on which the loop is based) and an arc $\gamma$ in $P \cap Q$. Suppose first that, $q$ lies in the annulus lying between some $\aaa_i$ and $\aaa_{i+1}$, $1 \leq i < m$.  $q$ can't be inessential in the annulus, since $|\bdd P \cap \bdd Q|$ has been minimized up to isotopy in $\bdd W$.  So $q$ spans the annulus.  The $1$--handle $W_i$  can be viewed as a regular neighborhood of  the arc $q$.  Then $D$ can be used to isotope $W_i$ through $\gamma \subset P$, removing both $\hat{\aaa}_i$ and $\hat{\aaa}_{i+1}$ and reducing $m$ by $2$.  A similar argument, reducing $m$ by $1$, applies if $q$ has one end on $\aaa_1$ and its other end on $p_0$.

The only remaining possibility is that $q$ is an essential arc in $\bdd W - \eta(\aaa)$ with both its ends on one of the components $\aaa_1$ or $\aaa_m$, say $\aaa_m$. (Both ends have to be on the same curve $\aaa_1$ or $\aaa_m$, since $\aaa$ is separating.)  That is, $q$ is an essential arc on the punctured torus $\bdd W_m - \hat{\aaa}_m$.  

\medskip

{\bf Claim}\qua  {\sl $q$ passes through a meridian of $W_m$ exactly once.}

\medskip
{\bf Proof of Claim}\qua The arc $\gamma$ together with a subarc $p$ of $\aaa_m$ form a closed curve in $P$, and so a closed curve that bounds a disk in $M[\aaa]$. If this disk is attached to $D$ along $\gamma$ the result is a disk $D_+$ in $M[\aaa]$  whose boundary $p \cup q$ is an essential closed curve on the boundary of $W_m$.  Since $M$ contains no Lens space summands or non-separating spheres, $p \cup q$ is a longitude of $W_m$, and so intersects a meridian of $W_m$ in a single point.  By isotoping $p$ to be very short, the intersection point lies in $q$, establishing the claim.

\medskip
Following the claim, $W_m$ can be viewed as a regular neighborhood of $q$, attached to $\hat{\aaa}_m$.  The disk $D_+$ appearing in the proof of the claim is a special $\bdd$--reducing disk for $W$.   The number of times $D_+$ passes through $\aaa$ is the number of vertices lying between $\gamma$ and $\aaa_m$ in $\hat{P}$ and so in particular is at most $m-1$.  This contradicts the original choice of $\hat{P}$.   The symmetric argument eliminates trivial loops in $\Sss$.  
\end{proof}

Following \fullref{lemma:noloops} we may immediately assume that $m, n \geq 2$ since, for example, if $m = 1$ the fact that fewer ends of edges of $\Sss$ lie on $p_0$ (if it exists) than on $\aaa_1$ guarantees that there will be a trivial loop at $\aaa_1$.  In the absence of trivial loops in the graphs, our analysis will focus on edges that are parallel in the graph.  Parallel edges cut off faces that are bigons. We will be interested in large families of parallel edges.  So consider a collection of parallel edges $e_1,\ldots,e_t$ in $\Uuu$, with ends on vertices $v, w$.  Number them in order around $v$, making an arbitrary choice between the two possible ways of doing this.  Arbitrarily call $v$ the source vertex.  Then the family of edges defines a function $\phi$ from a sequence of labels around $v$ (namely the labels of the ends of $e_1,\ldots,e_t$ at $v$) to a sequence of labels around $w$ (namely the labels of the ends of $e_1,\ldots,e_t$ at $w$).   In a typical setting we will know or assume a lot about the label sequence at $v$ (called the {\em source sequence}) and a little about the label sequence at $w$.  Typically we will only know that the label sequence at $w$ lies as a contiguous subsequence of a much larger sequence called the {\em target sequence}.  

In our setting, and momentarily assuming no labels $0$ appear, the labels around any vertex in $\Uuu$ appear in (circular) order $$1, \ldots, m, m, \ldots,1 ,1, \ldots, m, m, \ldots, 1 \hskip 3mm \ldots \hskip 3mm 1,\ldots,m,m,\ldots,1$$ 
with $|\bdd \aaa \cap \bdd \bbb | \geq 6$ determining the total number of sequences $1, \ldots, m$ and $m, \ldots ,1$ that appear.  When labels $0$ do appear, ie, when $\hat{P}$ is a special $\bdd$--reducing disk and $p_0$ intersects $\bdd \bbb$, then a label $0$ appears between some (but not necessarily all) successive labels $1, 1$ .  If the source sequence of a set of parallel edges (that is, the sequence of labels at the vertex $v$) is of length $t \leq m+1$ then the above long sequence is a natural target sequence.  That is, we know that the ends of the edges $e_1,\ldots,e_t$ at $w$ have labels some ordered contiguous subsequence in the long sequence above (with perhaps some $0$'s inserted).  But we could equally well have used the shorter target sequences $1,\ldots,m,m,\ldots,1,1,\ldots,m$ or $1,\ldots,m,m,\ldots,1,0,1,\ldots,m$ since any ordered sequence of $t \leq m+1$ contiguous labels in the long sequence is also an ordered contiguous sequence in $1,\ldots,m,m,\ldots,1,1,\ldots,m$ or $1,\ldots,m,m,\ldots,1,0,1,\ldots,m$.  

The following lemmas are classical, going back at least to Gordon and
Litherland \cite{GL} and Scharlemann \cite{Sc1}:

\begin{lemma}  \label{lemma:RP3}

Suppose $\phi$ is a function as described above determined by a set of parallel edges in $\Uuu$, with source sequence $0,1, \ldots, m$.  For all $0 \leq i \leq m-1$, if $\phi(i) = i+1$,   then $\phi(i+1) \neq i$.

\end{lemma}

\begin{proof}  The case $1 \leq i \leq m-1$ is representative.  The bigon between the edges represents a rectangle embedded in $M$ with one pair of opposite sides lying in $\hat{P}$. The other pair of opposite sides are parallel spanning arcs on the annulus in $\bdd W$ between $\aaa_i$ and $\aaa_{i+1}$.  Moreover these spanning arcs are oriented in such a way that if the rectangle between them in the annulus is added, the result is a M\"obius band $A$ with its boundary on $\hat{P}$.   The union of the M\"obius band $A$ and a disk component of $\hat{P} - \bdd A$ is a copy of $\mathbb{RP}^2$ in $M$, whose regular neighborhood is then a punctured $\Rpp$ in $M$, contradicting the hypothesis that $M$ contains no Lens space summands.  
\end{proof}  

\begin{lemma}  \label{lemma:incohere}

Suppose $\phi$ is a function as described above determined by a set of parallel edges  in $\Uuu$, with source sequence contained in $0,1, \ldots, m$.  Then:
\begin{itemize} 
\item For all $1 \leq i \leq m-1$, if $\phi(i) = i$,  then $\phi(i+1) \neq i-1$.
\item If $\phi(1) = 1$,   then $\phi(2) \neq 1$.  
\item If $\phi(m) = m$,   then $\phi(m-1) \neq m$. 
\end{itemize}
\end{lemma}

\begin{proof}  In each case, the hypothesis implies that an edge in $\Uuu$ with the same label on each end is incoherent.
\end{proof}

Jointly call these lemmas {\em the standard $\Rpp$ contradiction}, and note that they are typically applied to show that a set of mutually parallel edges can't be too numerous.  For example, in our context:

\begin{lemma}  \label{lemma:no1tom} No label sequence $0,1,\ldots,m$ or $1,\ldots,m, m$ appears as a source sequence for any set of parallel edges in $\Uuu$. Moreover, if no label $0$ appears in the target sequence, no label sequence $1,\ldots,m$ can appear as a source sequence.  
\end{lemma}

\begin{proof}   We prove the last statement first.  Consider the image of the label sequence $1,\ldots,m$ in the target sequence $1,\ldots,m,m,\ldots,1,1,\ldots,m$.   If $\phi(1)$  is the $\rrr^{th}$ term in the target sequence, call $0 \leq \rrr < 2m$, the {\em offset} of the function.   Consider the possible offsets.  

If the offset is trivial ($\rrr = 0$) or $\rrr = 2m$ then for all $1 \leq i \leq m$, $\phi(i) = i$.  This means that for each label $1 \leq i \leq m$ there is an edge in $\Uuu$ with that label at both ends.  Looking at the other graph, this means that every vertex in $\Sss$ is the base vertex for a loop.  An innermost loop then can have no vertices in its interior, ie, it would be a trivial loop, contradicting \fullref{lemma:noloops}. 

Suppose $1 \leq \rho \leq 2m-1$.  If $\rrr$ is even, then the label $m - \frac{\rrr}{2}$ contradicts \fullref{lemma:RP3}.  If $\rrr$ is odd then the label $m - \frac{\rrr-1}{2}$ contradicts \fullref{lemma:incohere}.  

The point is informally captured in the graph in \fullref{fig:no1tom}.

\begin{figure}[ht!]
\labellist\small
\pinlabel $m-\rho/2$ [r] at 55 137
\pinlabel $m$ [r] at 55 180
\pinlabel $\llap{\phantom{2}}m$ [t] at 186 48
\pinlabel $2m$ [t] at 316 48
\pinlabel $3m$ [t] at 447 48
\pinlabel \rotatebox{90}{label} [r] at 447 111 
\pinlabel {position in target sequence} [t] <0pt,-15pt> at 248 48
\endlabellist
\centering
\includegraphics[width=0.7\textwidth]{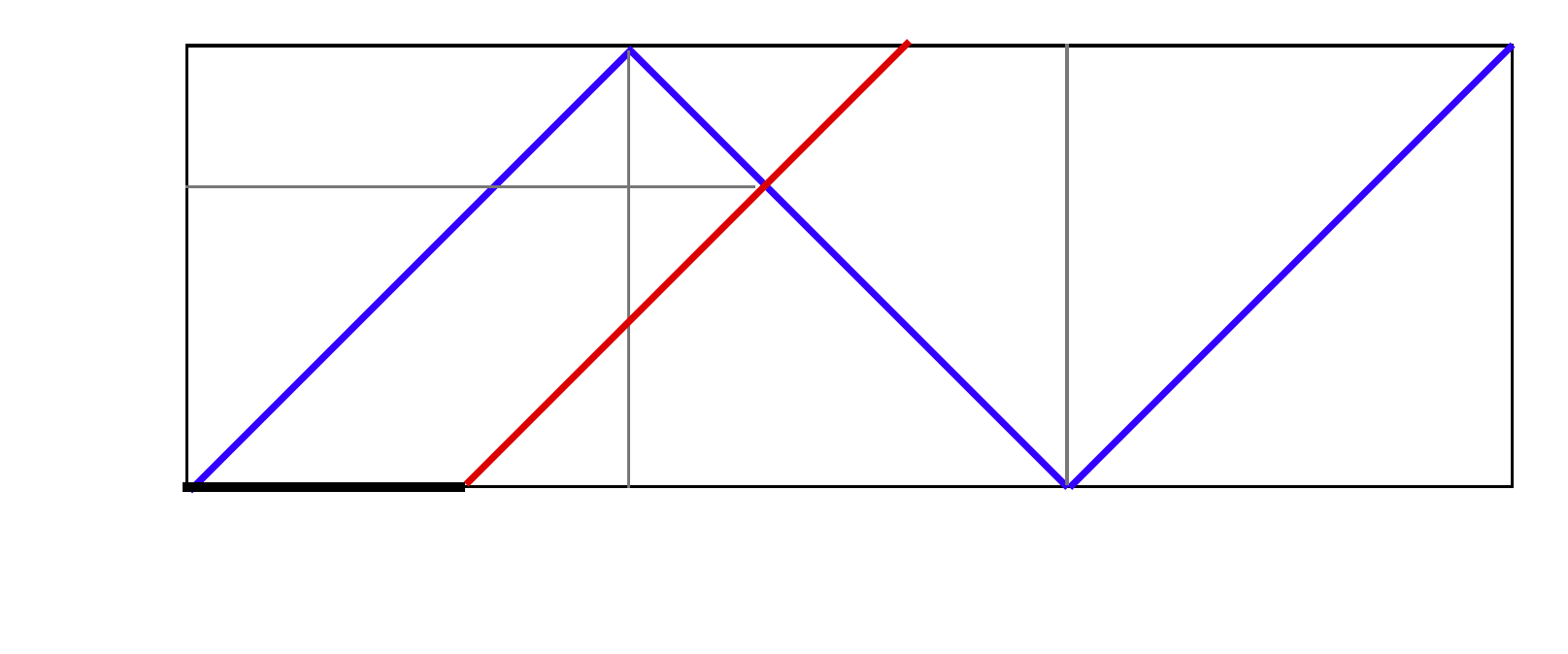}\vspace{-3mm}
\caption{} \label{fig:no1tom}
\end{figure}  

The case when the possible target sequence is  $1,\ldots,m,m,\ldots,1,0,1,\ldots,m$ is essentially the same argument, but requires using the entire source sequence $0, 1, \ldots, m$ or $1,\ldots,m,m$.  
\end{proof}

It is easy to complete the proof of \fullref{prop:bothsep} in the case that no labels $0$ appear (e. g. when $\hat{P}$ is a sphere or is a special $\bdd$--reducing disk whose boundary is disjoint from $\bbb$).  Since $|\bdd \aaa \cap \bdd \bbb | \geq 6$, every vertex in $\Uuu$ has on its boundaries at least six disjoint label sequences $1,\ldots,m$ or $m,\ldots,1$.  According to \fullref{lemma:no1tom} no such sequence can be entirely at the end of a parallel set of edges.  It follows that there are at least six breaks between sets of parallel edges at each vertex, hence at least six separate sets of parallel edges incident to each vertex.  This leads to a simple Euler characteristic contradiction, see \cite[Lemma 4.1]{GL} for the case in which $\hat{Q}$ is a sphere. 

The remaining case, when $\hat{P}$ is a special $\bdd$--reducing disk and $p_0$ intersects $\bdd \bbb$, is only slightly more delicate.  Note first that since $\bbb$ is separating, $|p_0 \cap  \bdd \bbb|$ is even and hence $\geq 2.$  This implies that each vertex of $\Uuu$ has at least two labels $0$ in its boundary.  Each label $0$ stands at the center of a sequence of edge labels $$ 1,\ldots,m,m,\ldots1,0,1,\ldots,m,m,\ldots,1$$ and each of these can be written as the end-point union of four label sequences  $$1,\ldots,m,m$$ $$m,\ldots,1,0$$ $$0,1,\ldots,m$$ $$m, m,\ldots,1.$$ Each of these four must have a break by \fullref{lemma:no1tom}.  Of the two sets of $4$ breaks thereby identified, at most two from each group can coincide since $|\bdd \aaa \cap \bdd \bbb | \geq 6$.  Hence there are again at least $8m - 2 = 6$ breaks around each vertex, leading to the same contradiction and
completeing the proof of the proposition.

\section{When meridians are non-separating -- an introduction}

Much complication is added if one or both of the meridians $\aaa$ or $\bbb$ is non-separating.  The most striking is that the graphs $\Sss$ and $\Uuu$ might have trivial loops, as we now describe.

Just as in the proof of \fullref{prop:bothsep}, no innermost trivial loop in $\Uuu$ can have ends labelled $i, i+1$, but there are three other possibilities:  as before they could have ends labelled $1, 1$ or $m, m$.  Or they could have ends labeled $1, m$ since, as $\aaa$ is non-separating, some arc of $\bdd \bbb - \bdd \aaa$ could have its ends on opposite sides of $\bdd \aaa$.  Only the last possibility (ends labeled $1, m$) withstands closer scrutiny:

\begin{lemma}  No innermost trivial loop in $\Uuu$ has both ends labelled $1$ (or, symmetrically, both ends labelled $m$).
\end{lemma}

\begin{proof}  Consider the twice punctured torus $T^* = \bdd W - \aaa$.  The two punctures (that is, $\bdd$--components) of $T^*$ can be identified with $\aaa_1$ and $\aaa_m$.  

\medskip
{\bf Claim 1}\qua  {\sl Any pair of arcs of $\bdd \bbb - \bdd \aaa$ that have both ends at $\aaa_1$ (or both ends at $\aaa_m$) are parallel in $T^*$.} 

\medskip
{\bf Proof of Claim 1}\qua Two non-parallel such arcs with ends at $\aaa_1$ would have complement in $T^*$ a punctured disk, with puncture $\aaa_m$.  All arcs of $\bdd \bbb - \bdd \aaa$ that have one end on $\aaa_m$ must then have their other end on $\aaa_1$ in order to be essential in $T^*$.  But then $\bdd \bbb$ would be incident to $\aaa_1$ at least four more times than it is incident to $\aaa_m$.  This is absurd, since $\aaa_1$ and $\aaa_m$ are parallel and $|\bdd P \cap \bdd Q|$ has been minimized up to isotopy.  The contradiction proves Claim 1.   

\medskip
{\bf Claim 2}\qua {\sl Any arc of $\bdd \bbb - \bdd \aaa$ that has
both ends at $\aaa_1$ or both ends at $\aaa_m$ is meridional in
$T^*$. That is, the arc, together with a subarc of $\aaa_1$ (or
$\aaa_m$) joining the ends, form a meridian circle on the boundary of
the solid torus $W - \aaa$.}

\medskip
{\bf Proof of Claim 2}\qua  Any outermost disk of $\bbb$ cut off by $\aaa$ intersects $T^* \subset \bdd W$ in an arc $\gamma$ which either has both ends at $\aaa_1$ or both ends at $\aaa_m$, say the former.  The outermost disk is a meridian of the solid torus $W - \aaa$ so $\gamma$, together with a subarc of $\bdd \aaa_1$ forms a meridian circle in the boundary of the solid torus $W - \aaa$.  This and Claim 1 establish Claim 2 for arcs with both ends at $\aaa_1$.

Any arc of $\bdd \bbb - \bdd \aaa$  with both ends at $\aaa_m$ (and a counting argument shows that there must be as many such arcs as there are with both ends at $\aaa_1$) must be meridional since if it had any other slope, it would necessarily intersect the meridional arc with both ends at $\aaa_1$.  This establishes Claim 2 also for arcs with both ends at $\aaa_m$.

\medskip

Following the claims, consider the disk $D$ cut off from $Q$ by an innermost trivial loop $\lambda$ in $\Uuu$ with both ends labeled $1$.  The two claims guarantee that the curve $\bdd D$ intersects $\bdd W - \aaa$ in a meridional arc.  In particular, the union of $D$ with the disk in $\hat{P}$ cut off by $\lambda$ and a subarc of $\bdd P$ is a disk $D_+$ in $M[\aaa]$ whose boundary in $W - \aaa$ is a meridian circle on the boundary of the solid torus $W - \aaa$.  The union of $D_+$ with a meridian disk for $W- \aaa$ contradicts the assumption that $M$ contains no non-separating spheres.  This completes the proof of the lemma.
\end{proof}

\begin{figure}[ht!]
\labellist\small
\pinlabel $\alpha_1$ [r] at 33 171
\pinlabel $\alpha_2$ [b] at 116 246
\pinlabel $\alpha_k$ [l] at 197 171
\pinlabel $\alpha_1^*$ [l] <0pt,-2pt> at 148 30
\pinlabel $\alpha_2^*$ [l] at 148 46
\pinlabel $\alpha_k^*$ [l] at 148 73
\endlabellist
\centering
\includegraphics[width=0.4\textwidth]{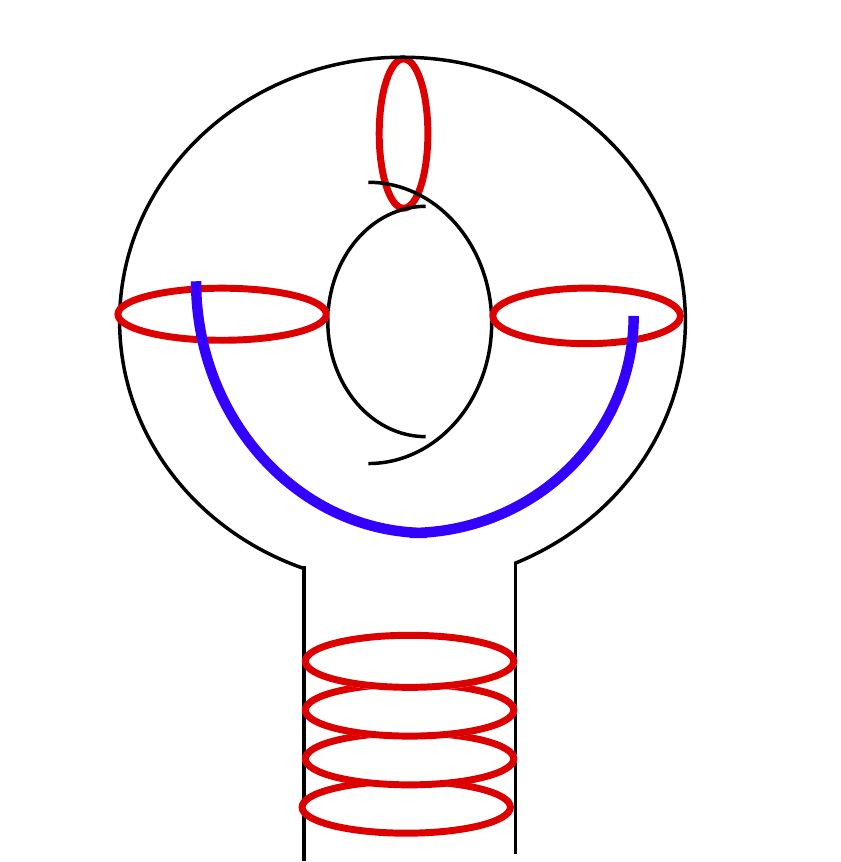}
\caption{} \label{fig:newvert}
\end{figure}

Following the lemma, we have that any innermost trivial loop in $\Uuu$
has ends labeled $1, m$.  The disk $D \subset Q$ it cuts off can be
used to $\bdd$--compress $P$ to $\bdd W$.  This alters $P$, replacing
the boundary curves $\aaa_1, \aaa_m$ with a separating curve bounding
a separating meridian $\aaa^{*}$ of $W$.  Replace $P$ then with a
planar surface (still called $P$) whose boundary components consist of
curves $\aaa_1,\ldots,\aaa_{k}$ parallel to $\bdd \aaa$ and curves
$\aaa^{*}_1,\ldots,\aaa^{*}_s$ parallel to $\bdd \aaa^{*}$, each
family labeled in order, and in such an order that one complementary
component in $\bdd W$ is a pair of pants with boundary components
$\aaa^{*}_s$, $\aaa_1$ and $\aaa_{k}$.  The construction guarantees
that $\bdd Q$ intersects the pair of pants in at least one arc with
ends on both $\aaa_1$ and $\aaa_{k}$, hence no arc with both ends on
$\aaa^{*}_s$.  See \mbox{\fullref{fig:newvert}.}

After this replacement, the vertices of $\Sss$ now are of two
different types and the labeling of edges around any vertex of $\Uuu$
is therefore more complicated.  But by allowing two vertex types,
choosing $P$ and $Q$ to minimize $|\bdd P \cap \bdd Q|$ now guarantees
that the corresponding graph $\Uuu$ in $\hat{Q}$ has no trivial loops.
Indeed, such a trivial loop would have its ends on one of:

\begin{itemize} 
\item $\aaa^*_1$, which would allow $s$ to be reduced by $1$, via
the argument of \fullref{lemma:noloops}, or
\item parallel copies of a meridian $\aaa$ or $\aaa^*$, which would
allow $s$ to be reduced by $2$, via the argument of
\fullref{lemma:noloops}, or
\item exactly two of the three meridians $\aaa_{k}, \aaa_1,
\aaa^*_s$.  In this case, a $\bdd$--compression removes the two
meridians incident to the loop and adds a copy of the third meridian,
reducing the total number of meridians and with it $|\bdd P \cap \bdd
Q|$.
\end{itemize}

The whole construction can be done symmetrically if $\bbb$ is non-separating:  Trivial loops could arise in $\Sss$ but be dealt with by identifying a separating meridian $\bbb^*$ and altering $Q$ so that it intersects $W$ in some curves parallel to $\bdd \bbb$ and some other curves parallel to $\bdd \bbb^{*}$.  

Although, after this alteration, the graphs $\Sss$ and $\Uuu$ have no trivial loops, the combinatorial argument now requires tracking four sets of circles.  This means that there are two sets of labels in each graph and two possible labeling schemes around vertices in each graph.  Also, the possible target sequences that arise are much more complicated than those that arise in the proof of \fullref{prop:bothsep}.  For example, a target sequence might contain $$1, \ldots, k, 1, \ldots, k, 1, \ldots, k, \ldots$$ intermixed with sequences of the form $$s^*, \ldots, 1^*, 1^*, \ldots, s^*.$$  As of this writing, such a combinatorial argument can be constructed for the case of reducing spheres. The argument appears likely to extend to the case of $\bdd$--reducing disks, but the addition of an extra boundary component (namely, the boundary component of the reducing disk) makes the final result (Conjecture 2 hence Conjecture 1) still uncertain.  

Scott Taylor \cite{Ta} has found a sutured manifold argument that avoids much of this combinatorial complication, in the same way that Gabai's \cite{Ga} circumvented the combinatorial difficulties of \cite{Sc3}.  From slightly different assumptions he is able to prove an analogue of Conjecture 2 that is even stronger (namely, the first conclusion holds unless $\alpha$ and $\beta$ are in fact disjoint) except in the most difficult but arguably the most interesting case: when $M = S^3$, $\alpha$ and $\beta$ are both non-separating and at least one of $M[\alpha]$ or $M[\beta]$ is an unknotted torus.

\bibliographystyle{gtart}
\bibliography{link}

\end{document}